\documentclass[final,onefignum,onetabnum]{siamart250211}

\usepackage{lipsum}
\usepackage{amsfonts}
\usepackage{graphicx} 
\usepackage{epstopdf}
\usepackage{algorithmic}
\usepackage{booktabs} %
\usepackage{amsopn}
\usepackage{tikz}
\usetikzlibrary{positioning}
\usepackage{standalone}
\usepackage{pgfplots}
\usepgfplotslibrary{groupplots}
\usetikzlibrary{patterns,patterns.meta}
\pgfplotsset{
	compat=1.18,
	scale only axis,
	every axis title shift=3pt,
	ymode=log,
	group style={group size=3 by 1, horizontal sep=0.75cm},
	height=1.5in,
	width=1.3in,
	xtick distance=10,
	minor x tick num=1,
	table/col sep=comma,
	table/search path={figs/},
	xlabel={$r$},
	title style={font=\small},
	label style={font=\footnotesize},
	tick label style={font=\footnotesize},
	legend style={font=\footnotesize},
	every axis plot/.append style={
		thick, 
		every mark/.append style={solid, fill opacity=0.3}
	},
	direct/.style={mark=star, black, dashed, mark size=3pt, every mark/.append style={ultra thick}},
	direct nonsym/.style={direct},
	direct sym/.style={mark=+, orange, dashed, mark size=3pt, every mark/.append style={ultra thick}},
    cg/.style={mark=+, cyan, dashed, mark size=3pt, every mark/.append style={ultra thick}, opacity=0.7},          %
	decoupled/.style={mark=square*, purple, dash dot},
	pcg/.style={mark=*, blue},
}

\usepackage{pgfplotstable}
\usepackage{booktabs,colortbl}
\usepackage{array}
\pgfplotsset{
	table/search path={.,figs},
}

\pgfplotstableset{
	error format/.style={
		column type={c},
		postproc cell content/.code={%
			\edef\doformat{\noexpand\pgfkeyssetvalue{/pgfplots/table/@cell content}{\noexpand\pgfmathprintnumber[fixed, zerofill, precision=3]{\pgfkeysvalueof{/pgfplots/table/@cell content}}}}%
			\doformat
		}
	},
	speedup format/.style={
		column type={@{}l},
		postproc cell content/.code={%
			\edef\doparens{\noexpand\pgfkeyssetvalue{/pgfplots/table/@cell content}{(\noexpand\pgfmathprintnumber[fixed,precision=0]{\pgfkeysvalueof{/pgfplots/table/@cell content}}$\noexpand\times$)}}%
			\doparens
		}
	},
	time format/.style={
		column type={>{~}r<{~~}},
		postproc cell content/.code={%
			\pgfplotstablegetelem{\pgfplotstablerow}{direct}\of{\timetable}
			\let\valdirect=\pgfplotsretval
			\pgfplotstablegetelem{\pgfplotstablerow}{direct_decoupled}\of{\timetable}
			\let\valdecoupled=\pgfplotsretval
			\pgfplotstablegetelem{\pgfplotstablerow}{pcg}\of{\timetable}
			\let\valpcg=\pgfplotsretval
			\pgfmathsetmacro{\minval}{min(\valdirect,\valdecoupled,\valpcg)}
			\pgfplotstablegetelem{\pgfplotstablerow}{\pgfplotstablecolname}\of{\timetable}
			\let\curval=\pgfplotsretval
			\pgfmathsetmacro{\ismin}{abs(\curval-\minval)<0.0001 ? 1 : 0}
			\ifnum\ismin=1
				\edef\doformat{\noexpand\pgfkeyssetvalue{/pgfplots/table/@cell content}{{\noexpand\color{blue}\noexpand\pgfmathprintnumber[fixed, zerofill, precision=2]{\pgfkeysvalueof{/pgfplots/table/@cell content}}}}}%
			\else
				\edef\doformat{\noexpand\pgfkeyssetvalue{/pgfplots/table/@cell content}{\noexpand\pgfmathprintnumber[fixed, zerofill, precision=2]{\pgfkeysvalueof{/pgfplots/table/@cell content}}}}%
			\fi
			\doformat
		}
	},
	align results/.style={
		columns={rank, err_direct, err_decoupled, err_pcg, direct, direct_decoupled, decoupled_speedup, pcg, pcg_speedup},
		every column/.style={font=\small, string type},
		columns/rank/.append style={column type={c|}},
		columns/err_direct/.append style={error format},
		columns/err_decoupled/.append style={error format},
		columns/err_pcg/.append style={error format, column type={r|}},
		columns/direct/.append style={time format},
		columns/direct_decoupled/.append style={time format},
		columns/decoupled_speedup/.append style={speedup format},
		columns/pcg/.append style={time format},
		columns/pcg_speedup/.append style={speedup format},
		every head row/.style={
			output empty row,
			before row={
				\toprule
				\multicolumn{1}{c|}{}& \multicolumn{3}{c|}{\textbf{Best Rel.~Error}} & \multicolumn{5}{c}{\textbf{Run Time (sec)}} \\
				\multicolumn{1}{c|}{$r$} & \multicolumn{1}{c}{Direct} & \multicolumn{1}{c}{Decoupled} & \multicolumn{1}{c|}{PCG}
				& \multicolumn{1}{c}{Direct} & \multicolumn{2}{c}{Decoupled} & \multicolumn{2}{c}{PCG} \\
				},
			after row=\midrule,
		},
		every last row/.style={after row=\bottomrule},
		every even row/.style={before row={\rowcolor[gray]{0.95}}},
	},
	time format unalign/.style={
		column type={>{~}r<{~~}},
		postproc cell content/.code={%
			\pgfplotstablegetelem{\pgfplotstablerow}{direct_sym}\of{\timetable}
			\let\valdirsym=\pgfplotsretval
			\pgfplotstablegetelem{\pgfplotstablerow}{direct_nonsym}\of{\timetable}
			\let\valdirnonsym=\pgfplotsretval
			\pgfplotstablegetelem{\pgfplotstablerow}{pcg}\of{\timetable}
			\let\valpcg=\pgfplotsretval
			\pgfmathsetmacro{\minval}{min(\valdirsym,\valdirnonsym,\valpcg)}
			\pgfplotstablegetelem{\pgfplotstablerow}{\pgfplotstablecolname}\of{\timetable}
			\let\curval=\pgfplotsretval
			\pgfmathsetmacro{\ismin}{abs(\curval-\minval)<0.0001 ? 1 : 0}
			\ifnum\ismin=1
				\edef\doformat{\noexpand\pgfkeyssetvalue{/pgfplots/table/@cell content}{{\noexpand\color{blue}\noexpand\pgfmathprintnumber[fixed, precision=0]{\pgfkeysvalueof{/pgfplots/table/@cell content}}}}}%
			\else
				\edef\doformat{\noexpand\pgfkeyssetvalue{/pgfplots/table/@cell content}{\noexpand\pgfmathprintnumber[fixed, precision=0]{\pgfkeysvalueof{/pgfplots/table/@cell content}}}}%
			\fi
			\doformat
		}
	},
	unalign results/.style={
		columns={rank, err_direct_sym, err_direct_nonsym, err_pcg, direct_sym, sym_speedup, direct_nonsym, nonsym_speedup, pcg},
		every column/.style={font=\small, string type},
		columns/rank/.append style={column type={c|}},
		columns/err_direct_sym/.append style={error format},
		columns/err_direct_nonsym/.append style={error format},
		columns/err_pcg/.append style={error format, column type={r|}},
		columns/direct_sym/.append style={time format unalign},
		columns/sym_speedup/.append style={speedup format},
		columns/direct_nonsym/.append style={time format unalign},
		columns/nonsym_speedup/.append style={speedup format},
		columns/pcg/.append style={time format unalign},
		every head row/.style={
			output empty row,
			before row={
				\toprule
				\multicolumn{1}{c|}{}& \multicolumn{3}{c|}{\textbf{Best Rel.~Error}} & \multicolumn{5}{c}{\textbf{Run Time (sec)}} \\[1mm]
				\multicolumn{1}{c|}{$r$} & \multicolumn{1}{c}{\shortstack{Direct\\Sym.}} & \multicolumn{1}{c}{\shortstack{Direct\\Nonsym.}} & \multicolumn{1}{c|}{PCG}
				& \multicolumn{2}{c}{\shortstack{Direct\\Sym.}} & \multicolumn{2}{c}{\shortstack{Direct\\Nonsym.}} & \multicolumn{1}{c}{PCG} \\
				},
			after row=\midrule,
		},
		every last row/.style={after row=\bottomrule},
		every even row/.style={before row={\rowcolor[gray]{0.95}}},
	},
}
\usetikzlibrary{matrix}
\usepackage{xcolor}
\usepackage{multirow}
\usepackage{subcaption,float}
\usepackage{tcolorbox} %
\tcbuselibrary{listings, skins, breakable} %

\usepackage{placeins} %
\usepackage{xurl} %

\usetikzlibrary{positioning} %

\usepackage{tensors}

\newcommand{\TheTitle}{Fast and Accurate CP-HIFI Tensor Decompositions Exploiting Kronecker Structure}
\newcommand{\TheShortTitle}{Fast CP-HIFI Tensor Decompositions}
\newcommand{\TheShortAuthors}{J.~J.~Brust and T.~G.~Kolda}

\headers{\TheShortTitle}{\TheShortAuthors}

\title{\TheTitle\thanks{Accepted for publication in SIAM J. Scientific Computing.\funding{JJB's work was partially supported by the startup fund at Arizona State University Grant PG16270 and the Simons Travel Support for Mathematicians.}}}

\author{
  Johannes J.~Brust%
  \thanks{
    School of Mathematical and Statistical Sciences, 
    Arizona State University, Tempe, AZ 
    (\email{jjbrust@asu.edu}).
  }
  \and 
  Tamara G.~Kolda%
  \thanks{
    MathSci.ai, Dublin, CA 
    (\email{tammy.kolda@mathsci.ai})
  }}

\DeclareMathOperator{\diag}{diag}

\renewcommand{\vec}{\operatorname{vec}}
\DeclareMathOperator{\unvec}{\operatorname{unvec}}

\NewDocumentCommand{\Real}{}{\mathbb{R}} 

\NewDocumentCommand{\bigO}{}{\mathcal{O}} 

\NewDocumentCommand{\A}{}{\Mx{A}} 

\NewDocumentCommand{\B}{}{\Mx{B}} 

\NewDocumentCommand{\Bbar}{}{\Mx[\bar]{B}} 

\NewDocumentCommand{\D}{}{\Mx{D}} 

\NewDocumentCommand{\Dbar}{}{\Mx[\bar]{D}} 

\NewDocumentCommand{\F}{}{\Mx{F}} 

\NewDocumentCommand{\G}{}{\Mx{G}} 

\NewDocumentCommand{\I}{}{\Mx{I}} 

\NewDocumentCommand{\Ihat}{}{\Mx[\hat]{I}} 

\NewDocumentCommand{\K}{}{\Mx{K}} 

\NewDocumentCommand{\Khat}{}{\Mx[\hat]{K}} 

\RenewDocumentCommand{\S}{}{\Mx{S}} 

\NewDocumentCommand{\Tobs}{}{\Tn{T}} 

\NewDocumentCommand{\Tk}{}{\Tm{T}{k}} 

\NewDocumentCommand{\T}{}{\Mx{T}} 

\NewDocumentCommand{\U}{}{\Mx{U}} %

\NewDocumentCommand{\V}{}{\Mx{V}} %

\NewDocumentCommand{\W}{}{\Mx{W}} 

\NewDocumentCommand{\Wbar}{}{\Mx[\bar]{W}} 

\NewDocumentCommand{\X}{}{\Mx{X}} 

\NewDocumentCommand{\x}{}{\Vc{x}} 

\NewDocumentCommand{\Z}{}{\Mx{Z}} 

\NewDocumentCommand{\Zhat}{}{\Mx[\hat]{Z}}

\NewDocumentCommand{\UK}{}{\Mx{U}{\Mx{K}}}
\NewDocumentCommand{\DK}{}{\Mx{D}{\Mx{K}}}
\NewDocumentCommand{\UV}{}{\Mx{U}{\Mx{V}}}
\NewDocumentCommand{\DV}{}{\Mx{D}{\Mx{V}}}

\ifpdf
\hypersetup{
  pdftitle={\TheTitle},
  pdfauthor={\TheShortAuthors}
}
\fi

\begin{document}

\maketitle

\begin{abstract}
Tensor decompositions are a fundamental tool in scientific computing and data analysis. In many applications---such as simulation data on irregular grids, surrogate modeling for parameterized PDEs, or spectroscopic measurements---the data has both discrete and continuous structure, and may only be observed at scattered sample points. 
The CP-HIFI (hybrid infinite-finite) decomposition generalizes the Canonical Polyadic (CP) tensor decomposition to settings where some factors are finite-dimensional vectors and others are functions drawn from infinite-dimensional spaces. 
The decomposition can be applied to a fully observed tensor (aligned) or, when only scattered observations are available, to a sparsely sampled tensor (unaligned). 
Current methods compute CP-HIFI factors by solving a sequence of dense linear systems arising from regularized least-squares problems to fit reproducing Kernel Hilbert space (RKHS) representations to the data, but these direct solves become computationally prohibitive as problem size grows. 
We propose new algorithms that achieve the same accuracy while being orders of magnitude faster. For aligned tensors, we exploit the Kronecker structure of the system to efficiently compute its eigendecomposition without ever forming the full system, reducing the solve to independent scalar equations. 
For unaligned tensors, we introduce a preconditioned conjugate gradient method, exploiting the problem's structure for fast matrix-vector products and efficient preconditioning. In our experiments, the proposed methods speed up the solution up to 500x compared to the prior naive direct methods, in line with the reduction in the theoretical computational complexity.
\end{abstract}

\begin{keywords}
reproducing kernel Hilbert space (RKHS), CP tensor decomposition, preconditioned conjugate gradients (pcg), Kronecker products
\end{keywords}

\begin{MSCcodes}
65F05,  %
15A69,  %
90C30,  %
\end{MSCcodes}

\section{Introduction}
\label{sec:intro}

The canonical polyadic (CP) tensor decomposition is a powerful tool for 
analyzing high-dimensional data \cite{CaCh70,Ha70,BaKo25}. 
It decomposes a tensor into a sum of outer products of finite-dimensional 
vectors, known as factors.
In recent work, Larsen et~al.~\cite{LaKoZhWi24} and 
Zhang et~al.~\cite{TaKoZh24} explore replacing some of the 
finite-dimensional factors (vectors) with infinite-dimensional factors (functions)
in a Reproducing Kernel Hilbert Space (RKHS) \cite{ScSm01}.
This approach is useful for data that can be viewed as sampled from a 
continuous distribution, such as time series or spatial data.
It is also useful for data that is not aligned on a regular grid, 
which may occur when measurements are irregular in time or come from
an adaptive spatial sampling process~\cite{TaKoZh24,LaKoZhWi24}.

Following \cite{LaKoZhWi24}, 
CP-HIFI (hybrid infinite- and finite-dimensional) tensor decomposition
is designed for a tensor that has one or more infinite-dimensional modes.
The data we observe is a finite set of samples,
represented as $\Tobs \in \Rmsiz$.
We let $\Omega \subseteq \mdom$ be the set of observations.
If we have all possible observations, 
then we say the problem is \emph{aligned}; 
otherwise, it is \emph{unaligned}.

The goal of CP-HIFI is to find a rank-$r$ CP approximation of the form
\begin{displaymath}
    \Tobs(\miwc) \approx \sum_{j=1}^r \A{1}(i_1,j) \, \A{2}(i_2,j)  \cdots  \A{d}(i_d,j),
\end{displaymath}
where $\A{k} \in \Real^{n_k \times r}$ is the factor matrix for mode $k$.
If mode $k$ is an infinite-dimensional mode, 
then we assume it is from an RKHS.
This means that its corresponding factor matrix is of the form
$\A{k} = \K{k} \W{k}$ where $\K{k} \in \Real^{n_k \times n_k}$ is a fixed positive semidefinite (psd) kernel matrix and $\W{k} \in \Real^{n_k \times r_k}$ is a
weight matrix.

There are a number of ways to solve this problem. Here, we focus on
solving for the infinite-dimensional modes in the 
alternating optimization approach as studied by 
Larsen et~al.~\cite{LaKoZhWi24}.
For an infinite-dimensional mode $k$, the subproblem to solve for $\W{k}$ is
\begin{equation}%
    \notag
    \min_{\W{k} \in \Real^{n_k \times r}} \nrm{ \Tk - \K{k}" \W{k}" \Z{k}'}_{\Omega}^2 + 
    \lambda_k \vec(\W{k})^{\Tr} (\I{r} \krn \K{k}) \vec(\W{k}),
\end{equation}
where $\Z{k}$ is the Khatri-Rao product of all factor matrices except the $k$-th one, $\lambda_k > 0$ is a regularization parameter, and $\nrm{\cdot}_{\Omega}$ is the norm over known entries.

\begin{figure}
    \centering
    \includegraphics{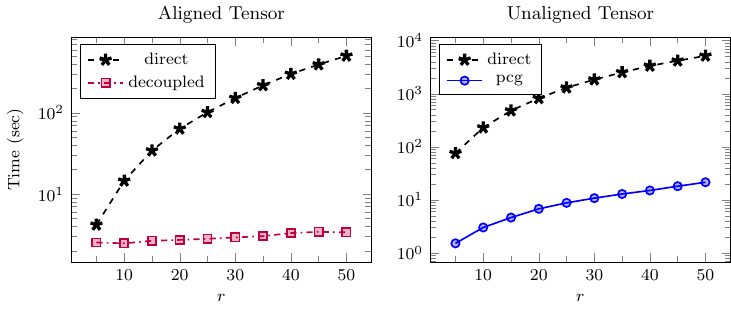}
    \caption{
        Comparison of direct methods proposed in prior work \cite{LaKoZhWi24} with the direct, decoupled, and pcg iterative methods proposed in this paper. The times represent the total time for an alternating optimization method. 
        See \cref{fig:miranda-aligned-comp,fig:miranda-unaligned} for further details.
    }
\label{fig:motivation}
\end{figure}

In \cite{LaKoZhWi24}, the authors solved the resulting linear systems of size $r n_k \times r n_k$ using direct methods. 
Setting $n \equiv n_k$, the direct method cost is $\bigO(n^3 r^3)$.
The goal of this paper is to investigate alternative solution methods.
In the aligned case, the problem can be decoupled, leading to an alternative direct solution for a solve cost of $\bigO(n^2r + nr^2 + r^3)$.
We also consider an iterative 
preconditioned conjugate gradient (PCG) \cite{GoVa13} solution to this problem,
for a cost of $\bigO(nr^2)$ per iteration.
In the unaligned case, pcg reduces the solve cost to $\bigO(rq + rn^2 + r^2n)$ per iteration, where $q$ is the number of known entries in $\Tobs$.
These provide orders of magnitude speedups over the naive direct methods,
as shown in \cref{fig:motivation}.
These methods also significantly reduce memory requirements.
We note that \cite{TaKoZh24} have proposed randomized solvers for these
subproblems, but these do not maintain the same accuracy. We discuss
options for randomized methods in the conclusions.

\paragraph{Notation}
Throughout this paper, we use the following notation.
Tensors are denoted by euler script letters (e.g., $\Tobs$),
matrices by bold uppercase letters (e.g., $\X$),
and vectors by bold lowercase letters (e.g., $\x$).
For a natural integer $n$, we define $[n] := \crly{1,2,\ldots,n}$.
A subscript on a matrix indicates an index into a sequence of matrices, e.g., $\A{k}$.
The $k$th mode unfolding of a tensor $\Tobs$ is denoted by a subscript in parentheses, e.g., $\Tk$.
Reviews of key concepts for tensors and Kronecker products are provided in \cref{app:tensor-background}, and more information can be found in \cite{BaKo25}.

\section{Review of CP-HIFI}
\label{sec:cp-hifi-review}

We briefly present the CP-HIFI method discussed in \cite{LaKoZhWi24}.
We assume that we have a finite-dimensional sample tensor 
\begin{displaymath}
	\Tobs \in \Rmsiz,
\end{displaymath}
where each $n_k$ corresponds to either the number of sample or design points (for infinite-dimensional modes) or simply the mode size (for finite-dimensional modes).

If every index is known, then we say the observations are \emph{aligned}; otherwise, we say they are \emph{unaligned}.
We can define the set of indices at which $\Tobs$ is known as
\begin{displaymath}
    \Omega = \crly*{ \Tobs(\miwc) \text{ is known} } \subseteq \mdom.
\end{displaymath}
We define 
$q \equiv |\Omega|$ to be the number of known entries
and the $\Omega$-norm to be 
$\nrm{\Tobs}_{\Omega}^2 \equiv \sum_{(i_1,i_2,\ldots,i_d) \in \Omega} \Tobs(i_1,i_2,\ldots,i_d)^2$.
In the aligned case, $\Omega = [n_1] \otimes [n_2] \otimes \cdots \otimes [n_d]$, $q=\prod_k n_k$, and $\nrm{\Tobs}_{\Omega} = \fnrm{\Tobs}$.

In CP-HIFI, we seek the factor matrices that minimize
\begin{equation}
    \nrm+{ \Tobs - \dsqr{\miwc[\A]}}_{\Omega}^2 + \sum_{k=1}^d \mathcal{R}_k,
\end{equation}
where 
the $\A{k}$ matrices affiliated with infinite-dimensional modes have a special structure
and
the $\mathcal{R}_k$ terms are regularization operators 
for the infinite-dimensional modes.
We solve this in an alternating fashion, optimizing for each factor matrix in  sequence while holding the others fixed.
The subproblem for mode $k$ is
\begin{equation}\label{eq:cp-hifi-subproblem}
    \min_{\A{k}} \nrm{ \Tk - \A{k}" \Z{k}'}_{\Omega}^2 + \mathcal{R}_k
    \qtext{where} \Z{k} := \skrp.
\end{equation}
Here, $\Tk \in \Real^{n_k \times M_k}$ is the mode-$k$ unfolding of $\Tobs$
and $M_k = \prod_{i \neq k} n_i$.

If $k$ is a finite-dimensional mode, then $\mathcal{R}_k=0$ and so \cref{eq:cp-hifi-subproblem} is equivalent to a CP-ALS minimization problem for which efficient solutions are known \cite{BaKo25}.

In this work, we focus on the infinite-dimensional subproblem. 
If $k$ is an infinite-dimensional mode, 
the factor matrix we are solving for is constrained 
to be of the form $\A{k} = \K{k}\W{k}$
where $\K{k} \in \Real^{n_k \times n_k}$ is a fixed psd kernel matrix
and $\W{k} \in \Real^{n_k \times r}$ is an unknown weight matrix.
The kernel matrix $\K{k}$ is computed by evaluating the bivariate kernel function (used to
define the RKHS) at all pairwise combinations of design points in mode $k$.
Solving for $\A{k}$ reduces to solving for $\W{k}$.
The regularization is of the form 
$\mathcal{R}_k = \vec(\W{k})^{\Tr} (\I{r} \krn \K{k}) \vec(\W{k})$
where $\I{r}$ is the $r \times r$ identity matrix. 
See \cite{LaKoZhWi24} for further details.
Thus, for an infinite-dimensional mode, \cref{eq:cp-hifi-subproblem} becomes 
\begin{equation}\label{eq:infinite-subproblem-with-k}
    \min_{\W{k} \in \Real^{n_k \times r}} \nrm{ \Tk - \K{k}" \W{k}" \Z{k}'}_{\Omega}^2 + 
    \lambda_k \vec(\W{k})^{\Tr} (\I{r} \krn \K{k}) \vec(\W{k}),
\end{equation}
where $\lambda_k > 0$ is a regularization parameter.

In the aligned case, all entries of $\Tobs$ are known, so
problem \cref{eq:infinite-subproblem-with-k} can be written as the regularized least-squares problem
\begin{equation}\label{eq:aligned-infinite-with-k}
    \min_{\W{k} \in \Real^{n_k \times r}} 
    \fnrm{ \Tk - \K{k} \W{k} \Z{k}'}^2 + \lambda_k \vec(\W{k})^{\Tr} (\I{r} \krn \K{k}) \vec(\W{k}).
\end{equation}

In the unaligned case, only some entries of $\Tobs$ are known. 
This is related to the problem of factorizing an incomplete tensor \cite[Ch.~14]{BaKo25}.
We let $\S{k}$ be the selection matrix of size $\prod_k n_k \times q$ such that $\S{k}'$ selects the $q$ known entries from the vectorization of the mode-$k$ unfolding of $\Tobs$. In particular, $\nrm{\Tobs}_{\Omega} = \tnrm{\S{k}' \vec(\Tk)}$.
Without loss of generality, we assume $\Tk$ is set to be zero at all unknown entries. 
In particular, this assumption means $\S{k}" \S{k}' \vec(\Tk) = \vec(\Tk)$.
Therefore, in the unaligned infinite-dimensional case, 
\cref{eq:infinite-subproblem-with-k} can be written as the regularized least-squares problem
\begin{equation}
    \label{eq:unaligned-infinite-with-k}
    \min_{\W{k} \in \Real^{n_k \times r}} 
    \tnrm{ \S{k}'(\Z{k} \krn \K{k})\vec(\W{k}) - \S{k}'\vec(\Tk) }^2 
    + \lambda_k \vec(\W{k})^{\Tr} (\I{r} \krn \K{k}) \vec(\W{k}).
\end{equation}

For notational convenience, we drop the subscript $k$ henceforth. 
The notation is as indicated in \cref{tab:notation}, including the sizes of each variable.
\begin{table}[ht]
    \centering\footnotesize
    \caption{Notation for mode-$k$ subproblem.}
    \label{tab:notation}
    \renewcommand{\arraystretch}{1.2}
    \pgfplotstabletypeset[
        col sep=&,row sep=\\,
        string type,
        columns/Symbol/.style={column type=c},
        columns/Description/.style={column type=l},
        columns/Size/.style={column type=c},
        every head row/.style={before row=\toprule, after row=\midrule},
        every last row/.style={after row=\bottomrule},
        every even row/.style={before row={\rowcolor[gray]{0.9}}},
        font=\footnotesize,
        assign column name/.style={/pgfplots/table/column name={\sffamily\textbf{#1}}},
    ]{
        Symbol & Description & Size \\
        $r$ & Target rank for CP-HIFI decomposition & scalar \\
        $\Omega$ & Set of known entries in $\mdom$ & set \\
        $q = |\Omega|$ & Number of known entries in $\Omega$ & scalar \\
        $\bar{n} = \sum_k n_k$ & Sum of sizes of all modes & scalar \\
        $n \equiv n_k$ & Size of mode $k$ & scalar \\
        $M \equiv \prod_{i \neq k} n_i$ & Product of sizes of all modes except $k$ & scalar \\
        $\bar{m} = \sum_{i \neq k} n_i$ & Sum of sizes of all modes except $k$ ($\bar{n} - n_k$) & scalar \\
        $\T \equiv \T{k}$ & Mode-$k$ unfolding of observed tensor, $\Tobs$ & $n \times M$ \\
        $\Z \equiv \Z{k}$ & Khatri-Rao product of all but the $k$th factor matrix & $M \times r$ \\
        $\K \equiv \K{k}$ & Kernel matrix for mode $k$ (infinite-dimensional) & $n \times n$ \\
        $\W \equiv \W{k}$ & Weight matrix for mode $k$ (infinite-dimensional) & $n \times r$ \\
        $\S \equiv \S{k}$ & Selection operator for $\Omega$ in mode $k$ (unaligned case) & $Mn \times q$ \\
    }
\end{table}

\section{Solving the Aligned Infinite-Dimensional (AI) CP-HIFI Subproblem}
\label{sec:aligned-subproblem}

In the aligned subproblem, dropping subscripts for ease of notation,
problem \cref{eq:aligned-infinite-with-k} becomes
\begin{equation}\label{eq:aligned-infinite}
    \min_{\W \in \Real^{n \times r}} 
    \fnrm{ \T - \K \W \Z' }^2 + \lambda \vec(\W)^{\Tr} (\I{r} \krn \K) \vec(\W).
\end{equation}
Following \cite{LaKoZhWi24}, 
if we calculate the gradient of \cref{eq:aligned-infinite} and set it to zero, we need to solve
\begin{equation}\label{eq:aligned-normal}
    \prn*{ (\Z'\Z \krn \K^2) + \lambda(\I{r} \krn \K) } \vec(\W) = \prn*{ \Z' \krn \K} \vec(\T).
\end{equation}
This system is symmetric positive semidefinite.
If we factor out $\I{r} \krn \K$ from \cref{eq:aligned-normal}, 
we arrive at a symmetric positive definite linear system of size $rn \times rn$:
\begin{equation}\label{eq:aligned-solve}
     \prn*{ \; \underbrace{\Z'\Z}_{\V} \, \krn \, \K  + \lambda \I{rn} }
     \vec(\W) = \vec(\; \underbrace{\T \Z}_{\B}\; ).  
\end{equation}
For ease of notation in the remainder of the section, we define
the MTTKRP
$\B := \T \Z \in \Real^{n \times r}$ and
the Gram matrix
$\V := \Z'\Z \in \Real^{r \times r}$.
Both $\B$ and $\V$ are needed for all methods discussed herein.
Computing $\B = \T \Z$ is an MTTKRP which costs $\bigO(Mnr)$ for dense $\T$,
and computing $\V = \Z'\Z$ is a Gram matrix that costs $\bigO(\bar{m} r^2)$ using the Khatri-Rao structure of $\Z$ \cite{BaKo25}.

We describe approaches for solving \cref{eq:aligned-solve}.
The direct method is described in \cref{sec:direct-ai} 
and is what was used in \cite{LaKoZhWi24}.
The direct solution can be greatly simplified by using 
the eigendecompositions of $\K$ and $\V$, leading to
a decoupled direct method in \cref{sec:direct-decoupled-ai}.
We also describe a small transformation of the system that allows for effective preconditioning and an iterative method in \cref{sec:pcg-ai}.
See \cref{tab:aligned-infinite-cost-comparison} and \cref{sec:cost-comparison-ai}
for a accounting of the costs and comparison of direct and iterative methods.

\subsection{Direct Solution of Original AI Subproblem}
\label{sec:direct-ai}
Since the matrix in \cref{eq:aligned-solve} is symmetric positive definite, 
we can use a direct method such as a Cholesky decomposition to solve it.
This is the method that was proposed in \cite{LaKoZhWi24}.

For the direct solve, we have to explicitly form the matrix
$\V \krn \K + \lambda \I $ of size $rn \times rn$, 
which costs $\bigO(r^2n^2)$
for the Kronecker product of two dense matrices.
The cost of solving the system via a Cholesky decomposition 
of the $rn \times rn$ matrix is $\bigO(r^3n^3)$.
The main memory cost is storing the matrix, which requires $\bigO(r^2 n^2)$ memory.

\subsection{Direct Decoupled Solution of AI Subproblem}
\label{sec:direct-decoupled-ai}
The problem can be solved much more efficiently using the eigendecompositions
of $\K$ and $\V$.

Let the eigendecomposition of the psd kernel matrix $\K$ be given by
\begin{equation}
    \label{eq:eigK}
    \K = \UK"  \DK"  \UK',
\end{equation}
where $\UK$ is an orthogonal matrix and $\DK$ is a diagonal matrix of nonnegative eigenvalues.
The factorization of $\K$ costs $ \bigO(n^3) $ flops.
This can be computed once and reused for all iterations of the alternating optimization method.

Let the eigendecomposition of the Gram matrix $\V$ be given by
\begin{equation}
    \label{eq:eigV}
    \V = \UV" \DV" \UV',
\end{equation}
where $\UV$ is an orthogonal matrix and $\DV$ is a diagonal matrix of nonnegative eigenvalues.
This will need to be computed once per iteration of the alternating optimization method, since $\V$ changes at each iteration. The cost is $\bigO(r^3)$.

Using \cref{eq:eigK,eq:eigV}, 
we can rewrite the system in \cref{eq:aligned-solve} as
\begin{equation}
    \prn(\UV \krn \UK) (\DV \krn \DK + \lambda \I{rn}) \prn(\UV \krn \UK)' \vec(\W) = \vec(\B).
\end{equation}
In other words, we have an eigendecomposition of the system.
This is a direct consequence of properties of Kronecker products: The eigendecomposition of a Kronecker product of two matrices is simply the Kronecker products of the eigendecompositions.
However, we can solve this without forming the Kronecker products explicitly.
The solution is
\begin{equation}\label{eq:W-decoupled}
    \W =  \UK \prn*{(\UK' \B \UV") \had \Mx{D}} \UV',
\end{equation}
where the matrix $\Mx{D}$ of size $n \times r$ is a reshape of the elementwise inverse of the diagonal of $\DV \krn \DK + \lambda \I{rn}$:
\begin{displaymath}
    \Mx{D} = \unvec( 1 / \diag(\DV \krn \DK + \lambda \I{rn}) ).
\end{displaymath}
Computing $\W$ via \cref{eq:W-decoupled} costs $\bigO(rn^2 + nr^2)$ flops,
corresponding to the matrix multiplies involving and $n \times r$ matrix with either $\UK$ or $\UV$. 
Additionally, all our computations are reduced to level-3 BLAS operations, which are highly optimized in modern linear algebra libraries.
The memory required is $\bigO(n^2 + r^2)$, corresponding to storing the eigendecompositions of $\K$ and $\V$.

\subsection{PCG Iterative Solution of Transformed AI Subproblem}
\label{sec:pcg-ai}  
As an alternative to direct solutions of the decoupled systems, 
we can use an iterative solver for \cref{eq:aligned-solve}.
We propose a transformation of the system
before we apply CG.

If we multiply \cref{eq:aligned-solve} 
on the left by the orthogonal matrix $\I{r} \krn \UK'$ and factor this matrix on the right,
we can rewrite \cref{eq:aligned-solve} as
\begin{equation}
    \prn{\V \krn \DK + \lambda \I{rn} } 
    \vec\prn{\, \underbrace{\UK'\W}_{\Wbar} \,}
    = \vec\prn{ \underbrace{\UK' \B}_{\Bbar} }.
\end{equation}
where we define $\Wbar := \UK'\W \in \Real^{n \times r}$ and 
$\Bbar := \UK'\B \in \Real^{n \times r}$.
The matrix-vector products with an arbitrary vector 
$\x$ with $\X = \unvec(\x) \in \Real^{r \times n}$
can be computed as
\begin{equation}\label{eq:pcg-matvec-ai}
    \prn*{\V \krn \DK + \lambda \I{rn}} \x
    = \vec\prn{\DK \X \V + \lambda \X}.
\end{equation}
Since $\DK$ is diagonal, each matrix-vector product
costs only $\bigO(nr^2)$.
Additionally, the computations are again reduced to highly-efficient level-3 BLAS operations.
We propose a diagonal preconditioner:
\begin{equation}
    \label{eq:precond-ai}
    \Dbar = \diag(\diag(\V)) \krn \DK  + \lambda \I.
\end{equation}
In other words, the diagonal matrix $\diag(\diag(\V))$ 
is constructed from the diagonal of $\V$.
The one-time cost to form $\Dbar^{-1}$ is $\bigO(nr)$, 
and the per-iteration cost to apply is also $\bigO(nr)$.
Once we have $\Wbar$, we can recover $\W$ via
$\W = \UK \Wbar$ for a cost of $\bigO(nr^2)$.

In our experiments, the iterative approach 
provides no advantage over the direct decoupled solution; however,
it serves as a basis for a similar approach for the unaligned case 
in \cref{sec:pcg-ui}.

\subsection{Comparison of Costs}
\label{sec:cost-comparison-ai}
A comparison of the direct, direct decoupled, and PCG iterative methods
are shown in \cref{tab:aligned-infinite-cost-comparison}.
For PCG, we let $p$ denote the number of iterations needed for convergence.
Recall that $n$ is the size of mode $k$, $r$ is the target rank,
$M$ is the product of the sizes of all modes except $k$,
and $\bar m$ is the sum of the sizes of all modes except $k$.
In general, we do not make assumptions about the relative sizes
of $n$, $r$, or $\bar m$. We do assume, however, that
$n,r,\bar{m} \ll M$.

\begin{table}[ht]
    \centering
    \caption{Comparison of costs to solve the mode-$k$ aligned infinite-dimensional subproblem \cref{eq:aligned-solve} of size $nr \times nr$ where $n$ is the size of mode $k$ and $r$ is the target tensor decomposition rank.
    The variables $M$ and $\bar m$ represent the product and sum of the sizes of all tensor modes except $k$, respectively.
    For the PCG iterative method, $p$ is the number of iterations.}
    \label{tab:aligned-infinite-cost-comparison}
    \renewcommand{\arraystretch}{1.5}
    \pgfplotstabletypeset[
        col sep=&,row sep=\\,
        string type,
        columns/Description/.style={column type=c<{\!\!\!}},
        columns/Direct/.style={column type=c, column name={Direct}},
        columns/DD/.style={column type=c, column name={Direct Decoupled}},
        columns/PCG/.style={column type=c, column name={PCG Iterative}},
        every head row/.style={before row={\toprule}, after row=\midrule},
        every last row/.style={after row=\bottomrule},
        every even row/.style={before row={\rowcolor[gray]{0.9}}},
        font=\footnotesize,
        assign column name/.style={/pgfplots/table/column name={\sffamily\textbf{#1}}},
    ]{
        Description & Direct & DD & PCG \\
        Factorize $\K$ \emph{one-time!} & --- & $\bigO(n^3)$ &  $\bigO(n^3)$ \\
        Compute MTTKRP $\B := \T \Z$ & $\bigO(rnM)$ & $\bigO(rnM)$ & $\bigO(rnM)$ \\
        Form Gram matrix $\V := \Z'\Z$ & $\bigO(r^2 \bar{m})$ & $\bigO(r^2 \bar{m})$ & $\bigO(r^2 \bar{m})$ \\
        Form $\V \krn \K + \lambda \I$ & $\bigO(r^2n^2)$ & --- & --- \\
        Factorize $\V$ & --- & $\bigO(r^3)$ & --- \\ 
        Solve system & $\bigO(r^3 n^3)$ & $\bigO(r^3 + r^2n + rn^2)$ & $\bigO( p r^2 n + rn^2)$ \\
        Total cost if $r < n$ & $\bigO(rnM + r^3n^3)$ & $\bigO(rnM + rn^2)$ & $\bigO(rnM + p r^2 n + rn^2)$ \\
        Storage & $\bigO(r^2 n^2)$ & $\bigO(r^2+n^2)$ & $\bigO(r^2+n^2)$ \\
    }
\end{table}

Overall, the direct solve of the original system has the highest computational cost as well as the highest memory requirement since it has to explicitly form and factorize a matrix of size $rn \times rn$.
This is the least efficient approach and should be avoided whenever possible.
The direct decoupled approach has the advantages of being very efficient to solve, having low memory requirements, avoiding forming large matrices, and having equivalent accuracy.
The PCG iterative has all the same advantages. The transformation of the problem and the simple preconditioner helps
it to achieve accuracy similar to the direct decoupled method.
Asymptotically, the direct decoupled should be more efficient than the PCG iterative method; however, these methods are similar in timing in the numerical experiments in \cref{sec:experiments}.
We have also included the iterative 
version here because it serves as a basis for the unaligned case in \cref{sec:pcg-ui}.

\section{Solving the Unaligned Infinite-Dimensional (UI) CP-HIFI Subproblem}
\label{sec:unaligned-subproblem}

In the unaligned case, only some entries of the tensor are known. 
Dropping subscripts for ease of notation,
problem \cref{eq:unaligned-infinite-with-k} becomes
\begin{equation}
    \notag
    \min_{\W \in \Real^{n \times r}} 
    \tnrm{ \S'(\Z \krn \K) \vec(\W) - \S' \vec(\T) }^2 
    + \lambda \vec(\W)^{\Tr} (\I{r} \krn \K) \vec(\W).
\end{equation}
Following \cite{LaKoZhWi24},
setting the gradient equal to zero results in the linear system
\begin{equation}\label{eq:unaligned-infinite} 
	\sqr*{\, 
        \underbrace{\prn(\Z \krn \K)' \S" }_{\F'}
        \underbrace{\S' \prn(\Z \krn \K)}_{\F} 
        + \lambda (\I{r} \krn \K) 
    } \vec(\W)
	= \prn(\I{r} \krn \K) \vec(\,\underbrace{\T\Z}_{\B}\,).
\end{equation}
We let $\B := \T \Z \in \Real^{n \times r}$ denote the MTTKRP, 
as we did in the aligned case.
We let $\F := \S' \prn(\Z \krn \K) \in \Real^{q \times nr}$
denote the subset of rows of $\Z \krn \K$ that correspond to the known entries of $\Tobs$.

We consider several approaches for solving \cref{eq:unaligned-infinite} 
in the remainder of this section.
We present a direct method for the symmetric linear system in \cref{sec:direct-ui-sym},
using an additional regularization term.
In \cref{sec:direct-ui-nonsym}, we present an alternative direct method
based on a nonsymmetric formulation of the system that does not
need additional regularization; this was the approach proposed in \cite{LaKoZhWi24}.
In \cref{sec:pcg-ui}, we present an iterative method based on
the symmetric system.
In \cref{tab:unaligned-infinite-cost-comparison} and \cref{sec:cost-compare-ui},
we provide an accounting of the costs and comparison of direct and iterative methods.

\subsection{Direct Solution of UI Subproblem}
\label{sec:direct-ui}

\Cref{eq:unaligned-infinite} is a 
symmetric positive semi-definite linear system of size $rn \times rn$.
Because of the structure of the unaligned problem, we cannot factor out
$(\I{r} \krn \K)$ as we did in the aligned case to arrive at a positive definite system.
Instead, since it is semi-definite, 
we add a regularization term parameterized by $\rho > 0$ 
to ensure positive definiteness.
The modified system is
\begin{equation}\label{eq:unaligned-infinite-1}
	\sqr*{ 
        \F' \F"  + \lambda (\I{r} \krn \K) + \rho \, \I{rn}
    } \vec(\W)
	= \vec(\K \B).
\end{equation}
Observe that we have pulled $\K$ inside the vectorization on the right-hand side.

\subsection{Symmetric Form}
\label{sec:direct-ui-sym}

To compute $\F$, we want to avoid forming the $Mn \times nr$ Kronecker product 
$\Z \krn \K$ explicitly.
Instead, we create two special matrices:
$\Khat \in \Real^{q \times n}$ and $\Zhat \in \Real^{q \times r}$.
Each index $\ell \in [q]$ corresponds to a known entry index that we 
denote as $(i_1^{(\ell)}, i_2^{(\ell)}, \dots, i_d^{(\ell)}) \in \Omega$.
Then, for each $\ell \in [q]$, we let
\begin{align}
    \label{eq:zhat} 
    \Zhat(\ell,:) & = \prn*{
        \A{d}(i_d^{(\ell)},:) \had \cdots \had \A{k+1}(i_{k+1}^{(\ell)},:) \had
        \A{k-1}(i_{k-1}^{(\ell)},:) \had \cdots \had \A{1}(i_1^{(\ell)},:)
    }, 
    \\
    \label{eq:khat} 
    \Khat(\ell,:) & = \K(i_k^{(\ell)},:)  .
\end{align}
In other words, $\Zhat$ and $\Khat$ 
represent the subset of rows of $\Z$ and $\K$, respectively,
that corresponds to the known entries of $\Tobs$.
Then, row $\ell$ of $\F$ is given by
\begin{equation}\label{eq:F}
    \F(\ell,:) = \Zhat(\ell,:) \krn \Khat(\ell,:).
\end{equation}
Making the typical assumption that the tensor order is less than any dimension, i.e., $d < n$, the computational complexity of forming $\F$ is $\bigO(qnr)$, and the cost of computing $\F'\F$ is $\bigO(qr^2n^2)$.

We omit the direct symmetric solver from our experimental comparisons in favor of the direct nonsymmetric solver described next in \cref{sec:direct-ui-nonsym}. 
There is no computational complexity advantage to the symmetric system over the nonsymmetric system.
Based on our testing, there is no computational advantage to solving the symmetric system over the nonsymmetric system for direct methods; to the contrary, the symmetric system is more expensive to solve because of the sensitivity to the additional regularization term $\rho\,\I{rn}$.

\subsection{Nonsymmetric Version}
\label{sec:direct-ui-nonsym}

An alternative direct solution method, which we include here because it was employed by \cite{LaKoZhWi24}, is to factor out $(\I{r} \krn \K)$ from \cref{eq:unaligned-infinite} to 
obtain the equation
\begin{equation} 
    \label{eq:unaligned-infinite-nonsym} 
	\sqr*{ \, \underbrace{ \prn(\Z \krn \I{n})' \S" }_{\G'} 
    \underbrace{\S' \prn(\Z \krn \K)}_{\F}  
    + \lambda \I{rn}} \vec(\W)
	=  \vec(\B),
\end{equation}
where we let $\G := \S' (\Z \krn \I{n}) \in \Real^{q \times rn}$
be analogous to $\F$.

We can form $\G$ similarly to how we formed $\F$.
If we define $\Ihat$ analogously to $\Khat$ so that $\Ihat(\ell,:) = \I{n}(i_k^{(\ell)},:)$ for each $\ell \in [q]$,
then row $\ell$ of $\G$ is given by
\begin{equation}\label{eq:G}
    \G(\ell,:) = \Zhat(\ell,:) \krn \Ihat(\ell,:).
\end{equation}
The computational complexity of forming $\G$ is $\bigO(qnr)$, and the cost of computing $\G'\F$ is $\bigO(qr^2n^2)$.

\subsection{PCG Solution of UI Subproblem}
\label{sec:pcg-ui}

As an alternative to the direct solvers, neither of which are particularly efficient for the unaligned case, we propose an iterative pcg solver for the symmetric system \cref{eq:unaligned-infinite-1}.

\subsubsection{Efficient Matrix-Vector Products for PCG Iterations}
We first consider efficient matrix-vector products, i.e., 
at every pcg iteration we need to compute the matrix-vector product:
\begin{align*}
    \Vc{y} &:= (\F' \F + \lambda (\I{r} \krn \K) + \rho \, \I{rn})\x \\
    &= \F' \F \x + \lambda (\I{r} \krn \K) \x + \rho \, \I{rn} \x \\
    &= \underbrace{\F' \F \x}_{\Vc{y}{1}} + \lambda \vec(\K\X) + \rho\x.
\end{align*}
Here $\x \in \Real^{rn}$ is an arbitrary vector, 
and $\X \in \Real^{n \times r}$ is such that $\vec(\X) = \x$.
The calculation of $\vec(\K\X)$ costs $\bigO(n^2r)$.
The matrix multiplication cost can potentially be reduced further if $\K$ has special structure; see, e.g., \cite{LiOnShCa20}. For example, $\K$ will be Toeplitz for a translation-invariant kernel on a uniform grid or circulant if the kernel is also periodic with period equal to the grid length.

We concentrate on computing $\Vc{y}{1} = \F' \F \x$ efficiently.
We will compute $\F' \F \x$ by first computing $\Vc[\bar]{x} := \F \x$ and then computing $\F'\Vc[\bar]{x}$. 
The matrix $\F$ is a subset of rows of the Kronecker product $\Z \krn \K$.
Conversely, $\F'$ is a subset of columns of the Kronecker product $\Z' \krn \K'$.
Pahikkala \cite{Pa14} and 
Airola and Pahikkala \cite{AiPa17} presented methods for computing matrix-vector products with matrices of this form, which we adapt here.%
\footnote{We became aware of this approach based on solutions recommended by AI systems as part of the ``First Proof'' challenge \cite{FP}. These solutions did not cite this prior work, but we found it after the fact. See the supplement for more information on the AI solutions and the connection to this prior work.}
Pahikkala \cite{Pa14} considers efficiently computing 
$\prn*{(\A \krn \B)\S\S'(\A \krn \B) + \lambda \I}\x$ where $\S$ is a selection matrix, which is exactly our problem.
Airola and Pahikkala \cite{AiPa17} consider the more general problem of efficiently computing $\S{1}'(\A \krn \B)\S{2}\Vc{x}$ where $\S{1}$ and $\S{2}$ are selection matrices. 

We break the computation into several steps, as follows:
\begin{align*}
    \Vc{y}{1} &= \F' \F \x = \prn*{\S' \prn(\Z \krn \K)}' \S' \prn(\Z \krn \K) \vec(\X)\\
    &= \prn*{\I{r} \krn \K} 
    \underbrace{\prn*{\Z' \krn \I{n}} \S 
        \underbrace{\S' \prn(\Z \krn \I{n}) 
            \underbrace{\prn(\I{r} \krn \K) \vec(\X)}_{\Vc[\hat]{x} = \vec(\K\X) \in \Real^{rn}}
        }_{\Vc[\bar]{x} \in \Real^{q}}
    }_{\Vc[\hat]{y}{1} \in \Real^{rn}}
    .    
\end{align*}
The computation of $\Vc[\hat]{x} = \vec(\K\X)$ is simply a matrix multiplication which costs $\bigO(n^2r)$.

The computation of $\Vc[\bar]{x} = \S' \prn(\Z \krn \I{n}) \Vc[\hat]{x}$ is more involved, but we can compute it efficiently by exploiting the structure of $\S$ and the Kronecker product per \cite{AiPa17}.
The $\ell$th row of $\S'$ picks out a single row of $\Z \krn \I{n}$,
specifically $\Zhat(\ell,:) \krn \Ihat(i_k^{(\ell)},:)$,  using the notation 
from \cref{eq:zhat,eq:khat}. Hence, for each $\ell \in [q]$, we have
\begin{equation}\label{eq:vl}
    v_{\ell} = \sqr*{\Zhat(\ell,:) \krn \I{n}(i_k^{(\ell)},:)} \vec(\Mx[\hat]{X})
    = \I{n}(i_k^{(\ell)},:) \Mx[\hat]{X} \Zhat(\ell,:)^{\Tr}
    = \Mx[\hat]{X}(i_k^{(\ell)},:) \Zhat(\ell,:)^{\Tr}.
\end{equation}
In other words, each entry of $\Vc[\bar]{x}$ is a dot product of a row of $\Zhat$ and a row of $\Mx[\hat]{X}$, so the cost to compute $\Vc[\bar]{x}$ is $\bigO(qr)$.

Computing $\Vc[\bar]{X}$ in two steps creates some efficiencies. If we were to instead compute $\Vc[\bar]{x} = \S' \prn(\Z \krn \K) \vec(\X)$ in a single step via $v_{\ell} = \K(i_k^{(\ell)},:) \Mx[\hat]{X} \Zhat(\ell,:)^{\Tr}$ for every $\ell \in [q]$, the total cost would be $\bigO(qnr)$.
Our proposed two-step computation costs $\bigO(qr + n^2r)$, which is more efficient under the natural assumption that $q > n$. A similar argument holds for the computation of $\Vc{y}{1}$ below.

With $\Vc[\bar]{x}$ computed, we can compute $\Vc[\hat]{y}{1} = \prn*{\Z' \krn \I{n}} \S \Vc[\bar]{x}$.
Define $\Mx[\bar]{X} \in \Real^{n \times M}$ such that $\vec(\Mx[\bar]{X}) = \S \Vc[\bar]{x}$, i.e., $\Mx[\bar]{X}$ is the same size as $\T$ and has zeros where $\T$ has missing entries.
Then, we are computing $\Vc[\hat]{y}{1} = \vec(\Mx[\bar]{X} \Z)$, which is really just a sparse MTTKRP if we think of $\Mx[\bar]{X}$ as being the mode-$k$ unfolding of a tensor
the same size and sparsity pattern as $\Tobs$. This can be computed efficiently by exploiting the structure of $\Z$ and the sparsity of $\Mx[\bar]{X}$.
Specifically, using $\Zhat$ from \cref{eq:zhat}, we have
\begin{equation}\label{eq:yhat1}
    \Mx[\hat]{Y}{1}(i,:) = \sum_{\ell: i_k^{(\ell)} = i} \bar{x}_{\ell} \, \Zhat(\ell,:),
\end{equation}
and then $\Vc[\hat]{y}{1} = \vec(\Mx[\hat]{Y}{1})$.
The cost to compute $\Vc[\hat]{y}{1}$ is $\bigO(qr)$, since we are summing over the $q$ known entries and each entry corresponds to a row of $\Zhat$ of size $r$.
Finally, we compute $\Vc{y}{1} = \prn*{\I{r} \krn \K} \Vc[\hat]{y}{1} = \vec(\K \Mx[\hat]{Y}{1})$, via a second matrix multiplication for a cost of $\bigO(rn^2)$.

Thus, the cost matrix-vector product cost per pcg iteration is $\bigO(n^2r + qr)$. 
In the next section, we discuss the preconditioner and its cost.

\subsubsection{Preconditioner for PCG Iterations}
\label{sec:precon-ui}
The main idea for the preconditioner is motivated by multiple similar
solutions produced by AI systems as part of the ``First Proof'' challenge \cite{FP}; full details are provided in the supplement.

The essential idea is to approximate the projection matrix $\S\S'$ by a multiple of the identity, i.e.,
$ \S \S' \approx \gamma \I $ where $ \gamma = q / N $ is the ``density'' of the unaligned tensor. 
Defining $\V := \Z'\Z$ as in the aligned case, the preconditioner 
is thus 
\begin{equation}\label{eq:precon}
    \Mx{M} 
    = \gamma \prn( \V \krn \K^2 ) 
    + \lambda \prn( \I{r} \krn \K ) 
    + \rho \I{rn}
\end{equation}
We can compute eigendecompositions of $\K$ and $\V$ as in \cref{eq:eigK,eq:eigV}.
This yields an orthogonal eigendecomposition of $\Mx{M}$:
\begin{equation}
    \Mx{M} = (\UV \krn \UK) \prn*{
        \gamma (\DV \krn \DK^2) + \lambda (\I{r} \krn \DK) + \rho \I{rn}
    } (\UV' \krn \UK').
\end{equation}
Define $\Vc{d} := 1 / \diag(\gamma (\DV \krn \DK^2) + \lambda (\I{r} \krn \DK) + \rho \I{rn}) \in \Real^{rn}$ as the diagonal of the \emph{inverse} of the diagonal middle matrix,
and let $\Mx{D} := \unvec(\Vc{d}) \in \Real^{n \times r}$.
Then, for an arbitrary vector $\x$, we can compute
\begin{equation}\label{eq:fx}
    \Mx{M}^{-1} \Vc{x} 
    = \UK" \prn{
    (\UK' \Mx{X} \UV") \had \Mx{D}
    } \UV',
\end{equation}
where $\had$ is the Hadamard (elementwise) product.

This preconditioner is a minor variant on the solution of the aligned system; 
see \cref{sec:direct-decoupled-ai}.
As in the aligned case, 
this requires a one-time eigendecomposition of $\K$ for a cost of $\bigO(n^3)$.
At each outer iteration, 
it requires an eigendecomposition of $\V$ for a cost of $\bigO(r^3)$.
The cost of applying the preconditioner per iteration is $\bigO(rn^2 + r^2n)$, corresponding to the matrix multiplies involving $\UK$ and $\UV$ in \cref{eq:fx}.
Additionally, all the preconditioner computations are reduced to highly-optimized level-3 BLAS operations.

\subsection{Comparison of Costs}
\label{sec:cost-compare-ui}

A comparison of 
the direct solution of the nonsymmetric problem \cref{eq:unaligned-infinite-nonsym},
and PCG iterative solutions of the 
symmetric problem \cref{eq:unaligned-infinite-1} 
are shown in \cref{tab:unaligned-infinite-cost-comparison}.
For PCG, we let $p$ denote the number of iterations needed for convergence.
Recall that $d$ is the order of the tensor, 
$n$ is the size of mode $k$, $r$ is the target rank,
and $q$ is the number of known entries.
In general, we do not make assumptions about the relative sizes
of $n$ and $r$. We do assume, however, that
$d < n,r \ll q$.
Because we are working with incomplete tensors, the MTTKRP is relatively cheap and never the dominant cost.

\begin{table}[ht]
    \centering
    \caption{Comparison of costs to solve the mode-$k$ unaligned infinite-dimensional subproblem \cref{eq:unaligned-infinite} of size $nr \times nr$ where $n$ is the size of mode $k$ and $r$ is the target tensor decomposition rank.
    The variable $q$ is the number of known entries in the observed tensor $\Tobs$.
    For the PCG iterative method, $p$ is the number of iterations.}
    \label{tab:unaligned-infinite-cost-comparison}
    \renewcommand{\arraystretch}{1.5} %
    \setlength{\tabcolsep}{4pt}
    \pgfplotstabletypeset[
        col sep=&,row sep=\\,
        string type,
        columns/Description/.style={column type=c},
        columns/DU/.style={column type=c, column name={Direct}},
        columns/PCG/.style={column type=c, column name={PCG Iterative}},
        every head row/.style={before row={
            \toprule
            }, after row=\midrule},
        every last row/.style={after row=\bottomrule},
        every even row/.style={before row={\rowcolor[gray]{0.9}}},
        font=\footnotesize,
        assign column name/.style={/pgfplots/table/column name={\sffamily\textbf{#1}}},
    ]{     
        Description & DU & PCG \\
        Factorize $\K = \U \D \U'$ \emph{one-time!} & --- &  $\bigO(n^3)$ \\
        Compute $\Zhat$ and MTTKRP $\B$ & $\bigO(drq)$ & $\bigO(drq)$ \\
        Form $\F$ and $\G$  & $\bigO(r n q)$ & --- \\
        Form matrix for linear solve &  $\bigO(r^2 n^2 q)$ & --- \\
        Form right-hand side & --- & $\bigO(r n^2)$  \\
        Form Preconditioner ($\Mx{M}$) & --- & $ \bigO(r^3 + rn)$ \\
        Solve system & $\bigO(r^3 n^3)$ &  $\bigO(pqr + pr^2n + prn^2)$ \\
        Total Cost if $d<r<n<q$ &  $\bigO(r^2 n^2 q + r^3 n^3)$ & $\bigO(prn^2+prq)$ \\ 
        Storage &  $\bigO(rnq + r^2 n^2)$ & $\bigO(rq)$ \\
    }
    \setlength{\tabcolsep}{6pt}
\end{table}

\paragraph{Summary and Comparison}

For the unaligned problem, there is no easy
decoupling as in the aligned case.
The direct method has cost which is at least cubic in the size of the unknown matrix $\W$.
In contrast, the PCG iterative method has a cost that is orders of magnitude lower, depending on the number of iterations $p$ needed for convergence and the relative sizes of $n$, $r$, $p$, and $q$. In general, we expect the problem to be well conditioned
so that $p$ is not too large.
The PCG method also has significantly lower storage requirements
than the direct methods.

\section{Numerical Experiments}
\label{sec:experiments}

We perform a series of numerical experiments to evaluate the performance of the new algorithms,
showing that the new methods are orders of magnitude faster than the direct methods proposed in \cite{LaKoZhWi24}.
We test our methods on two three-way datasets, continuous in all modes, described as follows
\begin{itemize}
    \item \textbf{Vortex tensor}: a $199\times 449\times 151$ tensor representation of a $199 \times 449 $ $x$-$y$ grid with 151 time snapshots \cite{KuBrBrPr16}. 
    This simulation demonstrates vortex shedding, a fluid dynamics phenomenon in which a flow separates behind an obstacle (e.g., a cylinder) and forms alternating vortices. 
    For the CP-HIFI factorizations, the RKHS uses Gaussian kernels with bandwidth parameters $\sigma_1=\sigma_2=4 $ and $ \sigma_3=3 $.
    \item
    \textbf{Miranda tensor}: a $ 256 \times 256 \times 269 $ tensor
    from the Scientific Data Reduction Benchmark (SDRBench) \cite{ZhDiLiLiTa20}
    showing one point in time in the 3-dimensional turbulent flow mixing of two fluids.
    For the CP-HIFI factorizations, the RKHS uses Gaussian kernels with bandwidth 
    parameters $ \sigma_1 = \sigma_2 = 1 $ and $ \sigma_3 = 2 $.
\end{itemize}
All numerical experiments are performed on a 64-bit Linux workstation with 128GB Ram, a
13th Gen Intel Core i9-13900KS CPU with 32 threads, and a NVIDIA GeForce RTX 4090 GPU. 
The tests are run on Matlab R2023a with the Tensor Toolbox (Version 3.8) \cite{BaKoTT25}.
Two repositories
contain the codes: 
\url{https://github.com/tgkolda/cp_hifi_code}
has implementations of the new algorithms and 
\url{https://github.com/tgkolda/cp_hifi_fast_exps}
has
drivers for experiments. 

\subsection{Numerical Experiments on Aligned Tensors}
\label{sec:aligned-experiments}

For aligned problems, we have observations for
every design point. For the example tensors, the design points are evenly spaced on a regular grid.
We compare three methods used to solve the subproblems
in the context of running CP-HIFI alternating optimization, as follows:
\begin{itemize}
    \item \textbf{direct}: direct solution adapted from \cite{LaKoZhWi24}, as described in \cref{sec:direct-ai}.
    \item \textbf{decoupled}: direct decoupled method from this paper in \cref{sec:direct-decoupled-ai}.
    \item \textbf{pcg}: iterative pcg method, also from this paper in \cref{sec:pcg-ai}. 
\end{itemize}
Each outer alternating optimization iteration 
updates all three factor matrices using the specified method.
Because CP-HIFI problems are nonconvex, each solver is run three times with different initial guesses for each rank $r \in \crly{5,10,\dots,50}$. 
We report results from the run with the lowest relative error.
In all experiments, we use the following parameters.
We use $\lambda=0.1$ as the RKHS regularization parameter in all modes for both datasets.
We set a limit of 50 maximum ALS iterations 
and a tolerance of $ 10^{-6} $ for the change
in the ALS error. 
For the iterative pcg solver, the maximum iterations are 75 and relative residual
tolerance is $ 10^{-6} $. 

\begin{figure}[ht] 
    \includegraphics{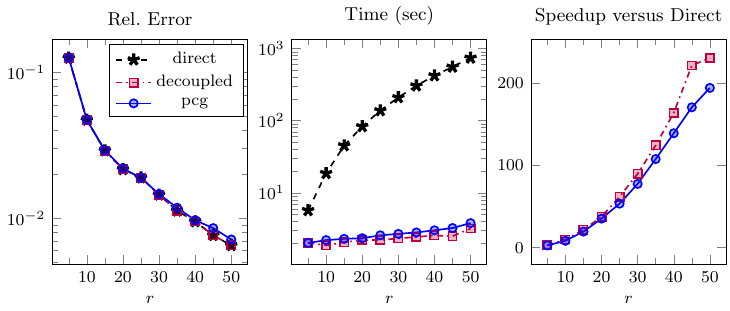}
    \caption{Aligned vortex tensor results for CP-HIFI, comparing relative error (left), runtime (middle), and speedup relative to the direct method (right)
    for different subproblem solvers.}
    \label{fig:vortex-aligned-comp}
\end{figure}

\begin{figure}[ht]
    \includegraphics{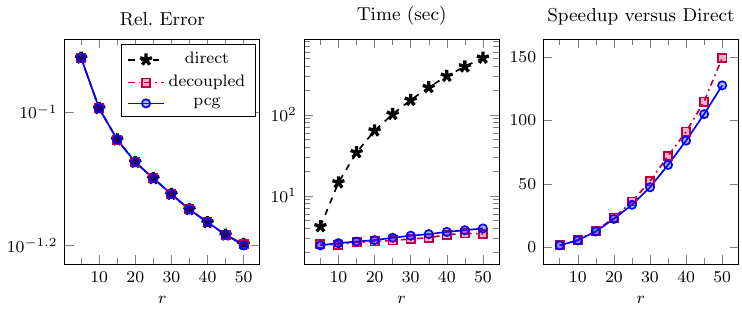}
    \caption{Aligned Miranda tensor results for CP-HIFI, comparing relative error (left), runtime (middle), and speedup relative to the direct method (right)
    for different subproblem solvers.}
    \label{fig:miranda-aligned-comp}
\end{figure}

Results for the aligned vortex and miranda tensors are shown in
\cref{fig:vortex-aligned-comp,fig:miranda-aligned-comp}, respectively. 
For both tensors and every choice of rank, the best solution ran for the full 50 
outer alternating optimization iterations. 
The left plots shows the relative error. 
The direct decoupled method achieves identical errors to the direct method;
the pcg method is only marginally worse.
The middle plots shows computational time. 
Both direct decoupled and pcg are significantly faster than the direct method. 
The right plots shows the speedup, which are as large as 200-fold for vortex and 130 fold for miranda.

\subsection{Numerical Experiments on Unaligned Tensors}
\label{sec:unaligned-experiments}
In the unaligned case, each tensor has only a subset of all possible observations. To simulate this with the datasets we have, we use a subsample of $q=50000$ entries from each tensor, randomly sampled from all possible entries. This equates to 0.37\% of the entries of the vortex tensor and 0.28\% of the entries of the Miranda tensor.
We compare the following methods which are used in the subproblem solver of CP-HIFI:
\begin{itemize}
    \item \textbf{direct}: As described in \cref{sec:direct-ui},
    we apply a standard direct method to a nonsymmetric formulation of the unaligned problem. This is the method used in \cite{LaKoZhWi24} and has the advantage of not requiring a second regularization parameter.
    \item \textbf{pcg}: Following logic similar to that of pcg for the aligned case, we develop a preconditioner, and apply pcg as described in \cref{sec:pcg-ui}. This is based on the symmetric version of the problem and uses an additional regularization term.
    \item \textbf{cg}: For comparison we also include cg without preconditioning.
\end{itemize}

As in the aligned case,
each outer alternating optimization iteration 
updates all three factor matrices using the specified method.
Also as in the aligned case,
because CP-HIFI problems are nonconvex, each solver is run three times with different initial guesses for each rank $r \in \crly{5,10,\dots,50}$. 
We report results from the run with the lowest relative error.
In all experiments, we use the following parameters.
We use $\lambda=0.1$ as the RKHS regularization parameter in all modes for both datasets.
For PCG we set $ \rho=10^{-6} $, which
did not cause any numerical problems in our experiments but yields accurate residuals.
We set a limit of 50 maximum ALS iterations 
and a tolerance of $ 10^{-6} $ for the change
in the ALS error. 
For the iterative solvers, the maximum ``inner'' iterations are 75 and relative residual
tolerance is $ 10^{-6} $. 

\begin{figure}[ht]
    \centering
    \includegraphics{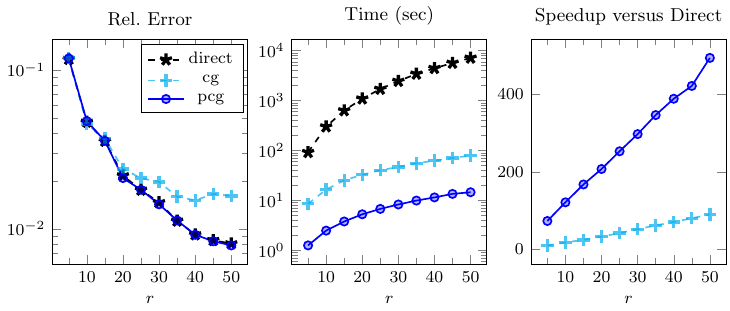}
    \caption{Unaligned vortex tensor results for CP-HIFI, comparing relative error (left), runtime (center), and the speedup of iterative methods versus the direct method.}
    \label{fig:vortex-unaligned}
\end{figure}

\begin{figure}[th]
    \centering
    \includegraphics{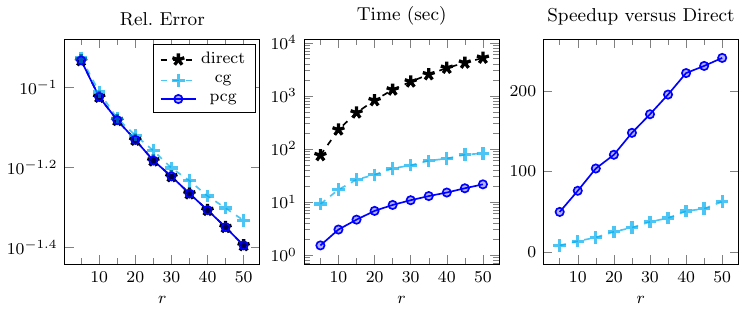}
    \caption{Unaligned Miranda tensor results for CP-HIFI, comparing relative error (left), runtime (center), and the speedup of the iterative solvers versus the direct method.}
    \label{fig:miranda-unaligned}
\end{figure}

\Cref{fig:vortex-unaligned,fig:miranda-unaligned} show results.
The leftmost plots show relative error, and the direct and pcg methods are generally comparable.
The cg method is less accurate for higher ranks.
The middle plots show runtimes. 
The right plots show speedups versus the direct method.
For the unaligned problems, the pcg speedup is even larger than for aligned problems. 
The cg method is significantly slower than pcg (due to taking more iterations to converge), but still faster than the direct method.

\begin{figure}[ht]
    \centering
    \graphicspath{{figs/}}
    \setkeys{Gin}{width=\textwidth}

    \pgfplotstableread[col sep=comma]{vortex_errs_unaligned.txt}\datatable
    \pgfplotstablegetelem{4}{pcg}\of\datatable
    \pgfmathprintnumberto[fixed,fixed zerofill,precision=4]{\pgfplotsretval}{\pcgerr}
    \pgfplotstablegetelem{4}{direct_nonsym}\of\datatable
    \pgfmathprintnumberto[fixed,fixed zerofill,precision=4]{\pgfplotsretval}{\dnserr}
    \pgfplotstableread[col sep=comma]{vortex_times_unaligned.txt}\datatable
    \pgfplotstablegetelem{4}{pcg}\of\datatable
    \edef\pcgtimeraw{\pgfplotsretval} %
    \pgfmathprintnumberto[fixed,precision=0]{\pcgtimeraw}{\pcgtime}
    \pgfplotstablegetelem{4}{direct_nonsym}\of\datatable
    \edef\dnstimeraw{\pgfplotsretval} %
    \pgfmathprintnumberto[fixed,precision=0]{\dnstimeraw}{\dnstime}
    \pgfmathparse{\dnstimeraw/\pcgtimeraw}
    \pgfmathprintnumberto[fixed,precision=0]{\pgfmathresult}{\speedup}

    \begin{tikzpicture}[
            node distance=10mm and 5mm, 
            picnode/.style={text width=2.00in,inner sep=0em},
            subcap/.style={below,text width=2.25in,inner ysep=0em},  
        ]
        \node[picnode] (orig) {\includegraphics[scale=0.95]{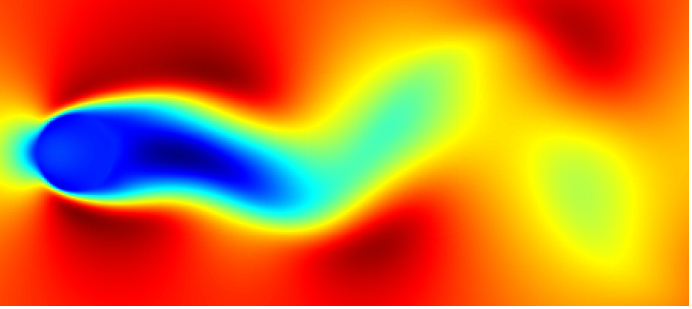}};
        \node[right=of orig,picnode] (samp) {\includegraphics[scale=0.95]{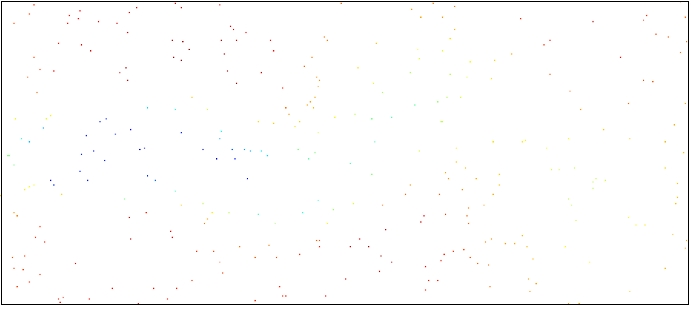}};
        \node[below=of orig,picnode] (pcg)  {\includegraphics[scale=0.95]{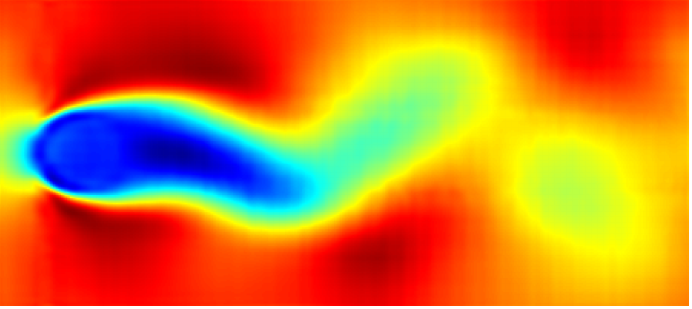}};
        \node[below=of samp,picnode] (nonsym) {\includegraphics[scale=0.95]{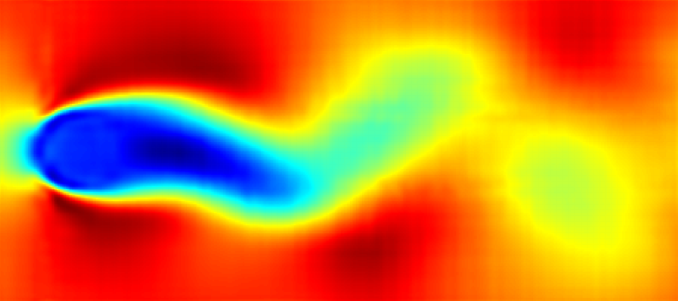}};
        \node[subcap] at (orig.south) {\subcaption{original}\label{fig:vortex-slices-orig}};
        \node[subcap] at (samp.south) {\subcaption{sampled}\label{fig:vortex-slices-samp}};
        \node[subcap] at (pcg.south) {\subcaption{pcg:\\error=\pcgerr, time=\pcgtime~sec}\label{fig:vortex-slices-pcg}};
        \node[subcap] at (nonsym.south) {\subcaption{direct:\\error=\dnserr, time=\dnstime~sec}\label{fig:vortex-slices-nonsym}};
    \end{tikzpicture}

    \caption{Panel (a) shows the $199 \times 449$ frontal slice 151 of the vortex tensor. Panel (b) shows the same slice of the sampled tensor (50,000 samples, i.e., 0.37\%). In (c)--(d) we compute CP-HIFI decompositions with rank $r=25$ using the subsampled data and different methods. The subsampled reconstructions are qualitatively similar, while pcg is more than \speedup~times faster than the direct solver.}
    \label{fig:vortex-slices}
\end{figure}

To compare the quality of the error differences, \cref{fig:vortex-slices} shows a reconstructed slice from the vortex tensor using the CP-HIFI model with $r=25$. \Cref{fig:vortex-slices-orig} shows the original unsampled slice 151, and \cref{fig:vortex-slices-samp} 
shows what remains in the sampled version of the tensor (50,000 observations in total, 0.37\% of the entries overall).
\Cref{fig:vortex-slices-pcg} shows the reconstruction based on the pcg solution, and \cref{fig:vortex-slices-nonsym}
shows the same for the direct nonsymmetric solution. These are qualitatively similar, but the pcg solution is significantly cheaper to compute.

\subsection{Numerical Experiments with Known Rank and Comparison to Missing Data Approach}
\label{sec:amino_acids}

As a final example, we consider the ability of CP-HIFI to recover the decomposition of a tensor of known rank.
In this case, we use the well-known amino acids fluorescence tensor from Bro \cite{Br97}.
This tensor is based on fluorescence spectrometer data from 201 emissions $\times$ 61 excitation wavelengths $\times$ 5 samples and is known to be rank $3$, which corresponds to the number of chemical components in the mixture.
The amino acids tensor is interesting for our purposes because the first two modes can be viewed as finite samples of a continuous process, so we can treat those modes as infinite-dimensional.
The standard CP-ALS algorithm applied to the full tensor yields a rank-3 decomposition corresponding to the chemical composition of the samples, and this will be used as ground truth. 

We attempt to recover the rank-3 decomposition with only \emph{limited observations} of the original tensor. 
For CP-HIFI, we treat the first two modes as infinite-dimensional and use Gaussian kernels to model the smoothness of the underlying functions.
For the first mode of dimension 201, we have $x$-values ranging from 1 to 201, and we use a Gaussian kernel with bandwidth $\sigma_1 = 15$.
For the second mode of dimension 61, we have $x$-values ranging from 1 to 61, and we use a Gaussian kernel with bandwidth $\sigma_2 = 8$.
We set the regularization parameter $\lambda = 0.1$ for both modes.
For comparison, we also consider a missing data approach without any continuity assumptions on the same limited observations.
There are a vast array of methods for missing data, but we choose a representative approach that is widely used and available in the Tensor Toolbox \cite{BaKoTT25}:
the weighted CP optimization method (CP-WOPT) developed by Acar et al.~\cite{AcDuKoMo11}, 
which uses an optimization approach to fit only the known entries of the tensor. 
We impose nonnegativity constraints on the factors because experiments without these constraints struggled with less than $q=5,000$ samples. 
\renewcommand{\floatpagefraction}{0.7}
\begin{figure}[pth]
    \centering
    \includegraphics{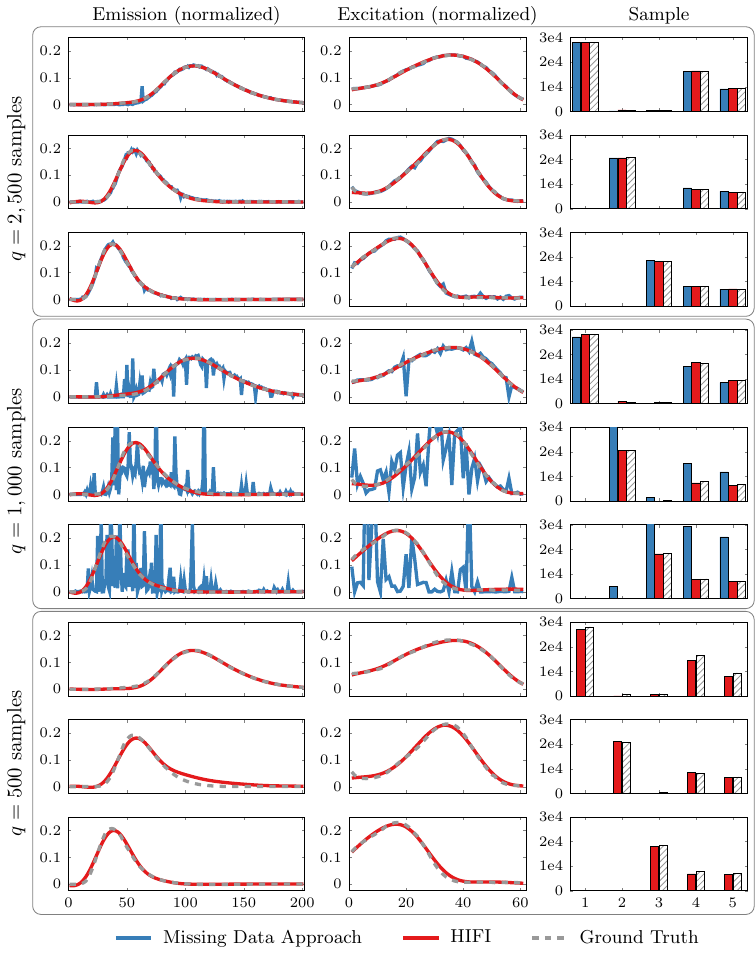}
    \caption{Comparison of a missing data approach and CP-HiFi on a tensor of size $201 \times 61 \times 5$ (61,305 total entries) with 2 modes that can be treated as infinite dimensional. We consider three different levels of missing data: $q=2,500$ (4\%), $q=1,000$ (1.6\%), and $q=500$ (0.8\%). The ground truth solution computed using the full dataset is shown as well. For $q=500$ the Missing Data Approach fails entirely and is not shown.}
    \label{fig:amino-acids}
\end{figure}

\Cref{fig:amino-acids} shows the results for CP-HIFI and CP-WOPT for three different sampling levels.
\begin{itemize}
    \item $\mathbf{q=2,500}$: Both CP-HIFI and CP-WOPT recover the rank-3 decomposition.
    The CP-WOPT solution has a few wiggles in its solution but is still easily interpretable.
    \item $\mathbf{q=1,000}$: CP-HIFI easily finds the ground truth solution while CP-WOPT struggles.
    CP-WOPT's first factor is fairly close to ground truth, but the wiggles in the solution are much more extreme than in the prior case. The second and third factors are not close to the ground truth.
    The problem is that the number of samples along each tensor fiber is too few, and there is no continuity ``glue'' to connect samples across fibers.
    \item $\mathbf{q=500}$: 
    CP-HIFI can still recover the rank-3 decomposition, though we see that it strays a bit from the true solution, especially in the second factor for the first mode.
    CP-WOPT fails entirely and is not even shown.
\end{itemize}

\par
We wouldn't generally expect that a missing data approach would work well with so few samples, so it is no surprise nor failing of the missing data approach per se. Instead, this example illustrates the power of CP-HIFI to work with exceptionally scarce data when an assumption of continuity can be used to weave together the limited information from different data points.

\pagebreak
\section{Conclusion}
\label{sec:concl}

This manuscript develops efficient solvers for the factors in a CP-HIFI decomposition. By exploiting a symmetric factorization
of the fixed kernel matrix and the gram of the Khatri-Rao matrices,
we can realize a significant reduction in computational cost as compared
to a naive approach. 
For unaligned tensors, we propose an efficient pcg iteration. 
The method uses the special structure of the problem to achieve efficient
matrix-vector products.
The inexpensive preconditioner is akin to the solution of the aligned case.
The unaligned pcg approach is based on solutions produced by AI systems as part of the ``First Proof'' challenge \cite{FP} and led us to prior work by Pahikkala \cite{Pa14} and Airola and Pahikkala \cite{AiPa17} on efficient matrix-vector products with Kronecker products and selection matrices.
Complexity analysis shows that the new methods scale orders of magnitude better
than the original methods, which is also born out in the numerical experiments. 
Our iterative and decoupled solvers may be extended in
future work. 

One topic for future work is handling nonnegativity constraints. 
For nonnegative tensors such as video data, 
it is typically desirable to constrain the factors to be nonnegative.
For CP-HIFI factors with nonnegative kernels, this reduces to computing nonnegative weights. 
However, the current approaches do not preserve positivity. Therefore, we need to investigate methods that are better suited for maintaining positivity; see, e.g., \cite{CoFaCo15}.
Another topic for future work is an ``all-at-once'' approach such as a Gauss-Newton method instead of alternating minimization.

For very large modes (i.e., $n$ is very large), factorizing the kernel matrix
may be prohibitively expensive. 
However, there are randomized method that can potentially accelerate that computation; 
see, e.g., \cite{WiSe00,FrTrUd23}.
Along the same line, our iterative methods may be combined with randomized sketching for applications to very large tensors \cite{BrSa24}.  
Randomized optimization methods have already been
applied to a related problem \cite{TaKoZh24}.
We generally assume that the sample points are given, but in some applications it may be possible to select the sample points and the choice of those sample points is another area for future research.

\appendix

\section{Background on Tensor Decompositions}
\label{app:tensor-background}
We refer the reader to \cite{BaKo25} for details of the operations discussed herein.
Let $\Tobs \in \Rmsiz$ be a $d$-way tensor. 
Its vectorization and mode-$k$ unfolding are denoted by
\begin{displaymath}
    \vec(\Tobs) \in \Real^{N}
    \qtext{and}
    \Tm{T}{k} \in \Real^{n_k \times M_k},
\end{displaymath}
respectively, where $N=\prod_k n_k$ and $M_k = N/n_k$.
The $\vec$ operation can also be applied to a matrix.

We let $\fnrm{\cdot}$ denote the matrix Frobenious norm and $\tnrm{\cdot}$ denote the vector 2-norm. The norm of a tensor is the square root of the sum of the squares of its entries:
\begin{displaymath}
    \nrm{\Tobs} = \tnrm{ \vec(\Tobs)} = \fnrm{ \Tm{T}{k}}.
\end{displaymath}

\subsection{Matrix Products}
\label{sec:matrix-products}
Certain matrix products will be used frequently, and their details can be found in \cite{BaKo25}.
The Hadamard product, denoted by $\had$, is the elementwise product.
The Kronecker product, denoted by $\krn$, is an outer product.
The Khatri-Rao product, denoted by $\krp$, is a columnwise Kronecker product.

\subsection{Kronecker Product Properties}
We use repeatedly Kronecker identities such as the following (c.f., \cite[Prop.~A.18]{BaKo25}):
\begin{align}
    \label{eq:kronecker-mult}
    \prn(\A \krn \B)' &= \A' \krn \B', \\
    (\A \krn \B)(\Mx{C} \krn \D) &= (\A\Mx{C}) \krn (\B\D), \\
    \label{eq:kronecker-identity}
    \vec(\A\B\Mx{C}) &= (\Mx{C}' \krn \A) \vec(\B)
\end{align}

Additionally, we use the Kronecker product perfect shuffle property, c.f., \cite[Prop.~A.19]{BaKo25}.
Let $\A \in \Real^{m \times p}$ and $\B \in \Real^{n \times q}$.
Then
\begin{equation}\label{eq:kronecker-shuffle}
    \A \krn \B = \Mx{P} (\B \krn \A) \Mx{Q},
\end{equation}
where $\Mx{P} \in \Real^{mn \times mn}$ and $\Mx{Q} \in \Real^{pq \times pq}$ are the $(m,n)$- and $(p,q)$-permutation matrices, respectively.
In particular, this means 
\begin{displaymath}
    (\A \krn \B) \vec(\X) = \vec(\Mx{Y})
    \quad \Leftrightarrow \quad
    (\B \krn \A) \vec(\X') = \vec(\Mx{Y}').
\end{displaymath}

\subsection{CP Decomposition and MTTKRP}
For a given $r \in \mathbb{N}$ and $\A{k} \in \Real^{n_k \times r}$ for $k \in [d]$,
we let $\dsqr{\miwc[\A]} \in \Rmsiz$ denote a rank-$r$ 
CP-decomposed tensor in Kruskal form \cite{BaKo25}. 
Its entries are given by
\begin{displaymath}
    \dsqr{\miwc[\A]}(\miwc) = \sum_{j=1}^r \prod_{k=1}^d \A{k} (i_k,j).
\end{displaymath}
Its mode-$k$ unfolding has special structure:
\begin{displaymath}
    \prn*{\dsqr{\miwc[\A]}}_{(k)} = \A{k} \prn*{\skrp}' \quad \in \Real^{n_k \times M_k}.
\end{displaymath}

The matricized-tensor times Khatri-Rao product (MTTKRP) \cite{BaKo25} is the product of 
the mode-$k$ unfolding of a tensor with the Khatri-Rao product of all
but the $k$th factor matrix:
\begin{displaymath}
    \Tm{T}{k} \prn*{\skrp} \quad \in \Real^{n_k \times r}.
\end{displaymath}

\section*{Acknowledgments}
A problem from this work was included in the First Proof Project \cite{FP},
and the AI-generated solutions to that problem have been incorporated into the methods described here. A description of our original solution and the 
solutions proposed by the AI, along with references to all the materials, is provided in the Supplemental Materials.

We thank the anonymous reviewers for their helpful comments and suggestions which have improved the quality of the manuscript.
 
\bibliographystyle{siamplain}


\end{document}